\newcommand{\calF}{\mathcal{F}}

\newcommand{\calO}{\mathcal{O}}

\newcommand{\calS}{\mathcal{S}}

\newcommand{\frakS}{\mathfrak{S}}

\newcommand{\bbF}{\mathbb{F}}

\newcommand{\bbP}{\mathbb{P}}
\newcommand{\bbQ}{\mathbb{Q}}

\newcommand{\bbZ}{\mathbb{Z}}

\newcommand{\xvar}{X_\varepsilon}
\newcommand{\PGL}{\textup{PGL}}
\newcommand{\PSL}{\textup{PSL}}
\newcommand{\SU}{\textup{SU}}

\newcommand{\Aut}{\textup{Aut}}
\newcommand{\Tr}{\textup{Tr}}
\newcommand{\rank}{\textup{rank}}
\newcommand{\Pic}{\textup{Pic}}

\def\Ker{{\text{Ker}}}

\def\PSU{{\textup{PSU}}}
\def\PGU{{\textup{PGU}}}
\def\GU{{\textup{GU}}}

\newcommand{\ord}{\text{ord}}

\newcommand{\Syl}{\text{Syl}}

\newcommand{\MW}{\text{MW}}

\documentclass[11pt]{amsart}
\usepackage{amscd, amsmath, amsthm, amssymb}
\newtheorem{theorem}{Theorem}[section]
\newtheorem{lemma}[theorem]{Lemma}
\newtheorem{proposition}[theorem]{Proposition}

\newtheorem{corollary}[theorem]{Corollary}
\theoremstyle{definition}     
\newtheorem{definition}[theorem]{Definition}

\theoremstyle{remark}
\newtheorem{remark}[theorem]{Remark}
\numberwithin{equation}{section}

\begin{document}
\title[K3 surfaces]
{K3 surfaces with a symplectic  automorphism of order 11}

\author[I. Dolgachev]{Igor V. Dolgachev}
\address{Department of Mathematics, University of Michigan, Ann Arbor, MI 48109, USA}
\email{idolga@umich.edu}
\thanks{Research of the first named author is partially supported by NSF grant
DMS-0245203}
\author[J. Keum]{JongHae Keum }
\address{School of Mathematics, Korea Institute for Advanced
Study, Seoul 130-722, Korea } \email{jhkeum@kias.re.kr}
\thanks{Research of the second named author is supported by KOSEF grant R01-2003-000-11634-0}
\begin{abstract}
We classify possible finite groups of symplectic automorphisms of
K3 surfaces of order divisible by 11.  The characteristic of the
ground field must be equal to 11. The complete list of such groups
consists of five groups: the cyclic group of order 11, $11\rtimes
5$, $L_2(11)$ and the Mathieu groups $M_{11}$, $M_{22}$. We also
show that a surface $X$ admitting an automorphism $g$ of order 11
admits a $g$-invariant elliptic fibration with the Jacobian
fibration isomorphic to one of explicitly given elliptic K3
surfaces.
\end{abstract}
\maketitle
\section{Introduction}

Let $X$ be a K3 surface over an algebraically closed field $k$ of
characteristic $p\ge 0$. An automorphism $g$ of $X$ is called
\emph{symplectic} if it preserves a regular 2-form of $X$. In
positive characteristic $p$, an automorphism of order a power of
$p$ is called \emph{wild}. A wild automorphism is symplectic. A
subgroup $G$ of the automorphism group $\Aut(X)$ is called
\emph{symplectic} if all elements of $G$ are symplectic, and
\emph{wild} if it contains a wild automorphism.

It is a well-known result of V. Nikulin that  the order of a
symplectic automorphism of finite order of a complex K3 surface takes value in the set
$\{1,2,3,4,5,6,7,8\}$. This
result is true over an algebraically closed field $k$ of positive
characteristic $p$ if the order is coprime to $p$. The latter condition is automatically satisfied if $p > 11$ \cite{DK2}.
If $p = 11$, a K3 surface $X_\varepsilon$ defined by the equation of degree 12 in $\bbP(1,1,4,6)$
\begin{equation}\label{formula}
y^2+x^3+\varepsilon x^2t_0^4+t_1^{11}t_0-t_0^{11}t_1 = 0, \quad
\varepsilon\in k,
\end{equation}
 admits a symplectic  automorphism of order 11
 \begin{equation}\label{form2}
 g_{\varepsilon}:(t_0,t_1,x,y) \mapsto (t_0,t_0+t_1,x,y).
 \end{equation}

The main result of the paper is the following.

\begin{theorem} \label{main} Let $G$ be a finite group of symplectic
automorphisms of a K3 surface $X$ over an algebraically closed
field of characteristic $p\ge 0$. Assume that the order of $G$ is
divisible by $11$. Then $p = 11$ and  $G$ is isomorphic to one of
the following five groups
$$C_{11},\, 11:5=11\rtimes 5,\, L_2(11) = \PSL_2(\bbF_{11}),\, M_{11},\, M_{22}.$$
Moreover, the following assertions are true.
\begin{itemize}
\item [(i)] For any element $g\in G$ of order $11$, $X$ admits a $(g)$-invariant
elliptic pencil $|F|$  and $X$ is $C_{11}$-equivariantly
isomorphic to a torsor of one of the surfaces $X_{\varepsilon}$
equipped with its standard elliptic fibration.
\item [(ii)]  If $X = X_\varepsilon$ and $G$ contains an element of order $11$
leaving invariant both the standard elliptic fibration and a
section, then $G\cong C_{11}$ if $\varepsilon\neq 0$ and
$G$ is isomorphic to a subgroup of $L_2(11)$ if $\varepsilon=0$.
\end{itemize}
\end{theorem}
The surface $X_0$ is a supersingular K3 surface with Artin
invariant 1 isomorphic to the Fermat surface
$$x_0^4+x_1^4+x_2^4+x_3^4 = 0.$$
In a recent paper of Kond\=o \cite{Ko} it is proven that both
$M_{11}$ and $M_{22}$ appear as symplectic automorphism groups of
$X_0$. An element $g$ of order $p=11$ in these groups leaves
invariant an elliptic pencil with no $g$-invariant section, and we
do not know whether the $g$-invariant elliptic pencil has no
sections or has a section but no $g$-invariant section. Thus the
surface $X_0$ admits three maximal finite simple symplectic groups
of automorphisms isomorphic to $L_2(11), M_{11}$ and $M_{22}$.

\begin{corollary} A finite group $G$ acts symplectically and
wildly on a K3 surface over an algebraically closed field of
characteristic $11$ if and only if $G$ is isomorphic to a subgroup
of $M_{23}$ of order divisible by $11$ and having $3$ or $4$
orbits in its natural action on a set of $24$ elements.
\end{corollary}

\bigskip {\it Acknowledgment}

\medskip
The authors are grateful to  S. Kond\=o for many fruitful discussions.

\bigskip
{\it Notation}

\medskip
For an automorphism group $G$ or an automorphism $g$ of $X$, we
denote by $X^g$ {\it the fixed locus} with reduced structure, i.e.
the set of fixed points of $g$.

A subset $T$ of $X$ is {\it $G$-invariant} if $g(T)=T$ for all
$g\in G$. In this case we say $G$ leaves $T$ invariant.

An elliptic pencil $|E|$ on $X$ is {\it $G$-invariant} if $g(E)\in |E|$ for
all $g\in G$. In this case we say $G$ leaves $|E|$ invariant.

\noindent We also use the following notations for groups:

\medskip
$C_n$ the cyclic group of order $n$, sometimes denoted by $n$,

$m:n=m\rtimes n$ the semi-direct product of cyclic groups $C_m$
and $C_n$,

$M_n$ the Mathieu group of degree $n$,

$\#G$ the cardinality of $G$,

$V^g$ the subspace of $g$-invariant vectors of $V$.

\section{The surfaces $X_0$ and $X_1$}

Let $p = 11$ and $X_\varepsilon$ be the K3-surface from
\eqref{formula}. The surface $X_\varepsilon$ has an elliptic
pencil defined by the projection to the $t_0,t_1$ coordinates
$$f_\varepsilon:X_\varepsilon\to \bbP^1.$$ We will refer to it as the {\it
standard elliptic fibration}. Its zero section, the section at
infinity, will be denoted by $S_\varepsilon$. It is immediately
checked that the surface $X_\varepsilon$ is nonsingular. Computing
the discriminant $\Delta_\varepsilon$ of the Weierstrass equation
of the general fibre of the elliptic fibration on $X_\varepsilon$
we find that
\begin{equation}\label{for2}
\Delta_\varepsilon =
-t_0^2(t_1^{11}-t_1t_0^{10})(5t_1^{11}-5t_1t_0^{10}+4\varepsilon^3
t_0^{11}).
\end{equation}
This shows that the set of singular fibres of the elliptic
fibration on $X_0$ (resp.  $X_{\varepsilon}, \varepsilon\neq 0$)
consists of 12 irreducible cuspidal curves (resp. one cuspidal
fibre and 22 nodal fibres). The automorphism $g_\varepsilon$ given
by \eqref{form2} is symplectic and of order $11$. It fixes
pointwisely the cuspidal fibres over the point $\infty = (0,1)$
and has 1 orbit (resp. 2 orbits) on the set of remaining singular
fibres. It leaves invariant the zero section $S_\varepsilon$. The
quotient surface $X_{\varepsilon}/(g_{\varepsilon})$ is a rational
elliptic surface with a double rational point of type $E_8$ equal
to the image of the singular point of the fixed fibre. A minimal
resolution of the surface has one reducible non-multiple fibre of
type $\tilde{E}_8$ and one irreducible singular cuspidal fibre
(resp. 2 nodal fibres).

\begin{proposition}\label{surface} Let $X$ be a K3 surface over an algebraically closed field $k$
of characteristic $11$. Assume that $X$ admits an automorphism $g$
of order $11$. Assume also that $X$ admits a $(g)$-invariant
elliptic fibration $f:X\to \bbP^1$ with a section $S$. Then there
exists an isomorphism $\phi:X\to X_\varepsilon$ of elliptic
surfaces such that $\phi g \phi^{-1} = \tau g_\varepsilon$ for
some translation automorphism $\tau$ of $X_\varepsilon$. In
particular, if $g(S)=S$ then $\phi g \phi^{-1}=g_\varepsilon$.
\end{proposition}

\begin{proof} Let $$y^2+x^3+A(t_0,t_1)x+B(t_0,t_1) = 0 $$ be the Weierstrass
equation of  the $g$-invariant elliptic pencil, where $A$ (resp.
B) is a binary form of degree  $8$ (resp. $12$). Since $f$ does
not admit a non-trivial $11$-torsion section (\cite{DK2},
Proposition 2.11), $g$  acts non-trivially on the base of the
fibration. After a linear change of the coordinates $(t_0,t_1)$ we
may assume that $g$ acts on the base by $$g:(t_0,t_1)\mapsto
(t_0,t_1+t_0).$$ We know that a $g$-invariant elliptic fibration
has one $g$-invariant irreducible cuspidal fibre $F_0$ and either
22 irreducible nodal fibres forming two orbits, or 11 irreducible
cuspidal fibres forming one orbit (\cite{DK1}, p.124). Thus the
discriminant polynomial $\Delta = -4A^3-27B^2$ must have one
double root (corresponding to the fibre $F_0$) and either one
orbit of double roots or two orbits of simple roots.  We know that
the zeros of $A$ correspond to either cuspidal fibres or
nonsingular fibres with ``complex multiplication'' automorphism of
order 6. Since this set is invariant with respect to our
automorphism of order 11 acting on the  base, we see that the only
possibility is $A = ct_0^8$ for some constant $c\in k$. We obtain
$\Delta = -4c^3t_0^{24}-27B^2.$ Again this uniquely determines $B$
and hence the surface. Since $B$ is of degree 12 and invariant
under the action of $g$ on the base, it must be of the form
$$B=a(t_1^{11}-t_1t_0^{10})t_0+bt_0^{12},$$
for some constants $a$, $b$. One can rewrite the above Weierstrass
equation in the form $$y^2+x^3+\varepsilon
x^2t_0^4+a(t_1^{11}t_0-t_0^{11}t_1)+b't_0^{12} = 0.$$ A suitable
linear change of variables $u_0=t_0,\, u_1=t_1+dt_0$ makes $b'=0$
without changing the action of $g$ on the base. Thus $X \cong
X_\varepsilon$ as an elliptic surface. Let $\phi:X\to
X_\varepsilon$ be the isomorphism. The composite
$$\phi g \phi^{-1}g_\varepsilon^{-1}:X_\varepsilon\to
X_\varepsilon$$ acts trivially on the base, hence must be a
translation automorphism. Since $\phi$ maps the zero section $S$
of $f:X\to \bbP^1$ to the zero section $S_\varepsilon$ of
$f_\varepsilon:X_\varepsilon\to \bbP^1$ and
$g_\varepsilon(S_\varepsilon)=(S_\varepsilon)$, the last assertion
follows.
\end{proof}

\begin{lemma}\label{taug} Let $\varepsilon=0$. For any translation automorphism $\tau$ of $X_0$,
the composite automorphisms $\tau g_0$ and $g_0\tau$ are of order
$11$.
\end{lemma}

\begin{proof} Let $f:X\to B$ be any elliptic surface with a section $S$. Recall
that its Mordell-Weil group $\MW(f)$ is isomorphic to the quotient
of the Neron-Severi group by the subgroup generated by the divisor
classes of $S$ and the components of fibres. Thus any automorphism
$g$ of $X$ which preserves the class of a fibre and the section
$S$ acts linearly on the group $\MW(f)$. Assume $\MW(f)$ is
torsion free.  Suppose $g$ is of finite order $n$ with
$\rank~\MW(f)^g=0$ and let $\tau$ be a translation automorphism
identified with an element of $\MW(f)$. Then, for any $s\in
\MW(f)$ we have
$$\tau g(s) = g(s)+\tau, \quad  (\tau g)^n(s) = g^n(s)+g^{n-1}(\tau)+\ldots +g(\tau)+\tau = s.$$
The last equality follows from that the linear action of $g-1_X$
on $\MW(f)$ is invertible. This shows that $(\tau g)^n$ acts
identically on $\MW(f)$. It also acts identically on the class of
a fibre. Thus $(\tau g)^n$ acts identically on the Neron-Severi
lattice.

Apply this to our case $\varepsilon=0$, when $g = g_0$ is a
symplectic automorphism of order 11 of $X_0$. We will see in the
proof of Proposition \ref{Hmax} that $\MW(f_0)$ is torsion free.
By Lemma \ref{lem1}(iii) below, $\rank~\MW(f_0)^g=0$. Since the
surface $X_0$ is supersingular (see Remark \ref{goto}), by a
theorem of Ogus \cite{Ogus}, an automorphism acting identically on
the Picard group must be the identity. Thus $\tau g_0$ is a
symplectic automorphism of order 11 for any section $\tau$.
\end{proof}

An interesting question: Is there a $\tau$ such that the fixed
locus $X_0^{\tau g_0}$ consists of an isolated point, the cusp of
a cuspidal curve fixed pointwisely by $g_0$? We do not know any
example of a symplectic automorphism of order 11 with an isolated
fixed point.

\begin{lemma}\label{lem1} Let $X$ be a K3 surface over an algebraically closed field $k$
of characteristic $11$. Assume that $X$ admits an automorphism $g$
of order $11$. Then the following assertions are true.
\begin{itemize}
\item [(i)] $X$ admits a $(g)$-invariant
 elliptic pencil $|F|$;
\item [(ii)] $\rank~\Pic(X/(g))=2$;
\item [(iii)] for any $l\ne 11$, $\dim H^2_{\rm et}(X,{\bbQ}_l)^g=\rank~\Pic(X)^g=2;$
\item [(iv)] $\rank~\Pic(X)=2$, $12$ or $22$.
\end{itemize}
\end{lemma}

\begin{proof} To prove (i), assume first that
$X$ does not admit a $(g)$-invariant elliptic pencil and $X^g$ is
a point. This case could happen only if the sublattice $N$ of the
Picard group of a minimal resolution of $X/(g)$ generated by
irreducible components of exceptional curves is $11$-elementary,
and $N^\perp$ is an even lattice of rank 2. This is contained in
the proof of Proposition 2.9 of \cite{DK2}. The intersection
matrix of $N^\perp$ is of the form
\[\begin{pmatrix}2a&c\\
c&2b\end{pmatrix}\]
 Since $N^\perp$ is indefinite and $11$-elementary,
$$\det N^\perp = 4ab-c^2=-1, \quad -11\quad {\rm or}\quad -121.$$
In the first case, $N^\perp\cong U$, where $U$ is an even
indefinite unimodular lattice. The second case cannot occur, since
no square of an integer is congruent to 3 modulo 4. Assume the
third case. Since $N^\perp$ is $11$-elementary, all of the
coefficients of the matrix must be divisible by 11, and hence
$N^\perp\cong U(11)$. Therefore, in any case $N^\perp$ contains an
isotropic vector. This is enough to deduce that $X$ admits a
$(g)$-invariant elliptic pencil by the same proof as in
Proposition 2.9 of \cite{DK2}.

Let $|F|$ be a $(g)$-invariant elliptic pencil. It follows from
\cite{DK1}, p. 124, that the elliptic fibration has one cuspidal
fibre and 22 nodal fibres, or 12 cuspidal fibres. The automorphism
$g$ leaves one cuspidal fibre $F_0$ over a point $s_0\in \bbP^1$
invariant.

Assertion (ii) follows from \cite{DK1}, where we proved that
$X/(g)$ is a rational elliptic surface with no reducible fibres,
and its minimal resolution is an extremal elliptic surface, i.e.
the sublattice of the Picard group generated by irreducible
components of fibres is of corank 1.

It is proven in \cite{HN}, Proposition 3.2.1, that for any $l\ne
p$ coprime with the order of $g$
$$\dim H^2_{\rm et}(X,{\bbQ}_l)^g=\dim H^2_{\rm et}(X/(g),{\bbQ}_l).$$
In fact it is true for all $l\ne p$ because of the invariance of
the characteristic polynomial of an endomorphism of a smooth
algebraic variety. Now by (ii),
$$\dim H^2_{\rm et}(X,{\bbQ}_l)^g=\dim H^2_{\rm
  et}(X/(g),{\bbQ}_l)=\rank~ \Pic(X/(g))=2.$$
Since $g$ fixes the class of a fibre and an ample divisor,
$\rank~\Pic(X)^g\ge 2.$ This proves (iii).

Considering the $\bbQ$-representation of the cyclic group $(g)$ of
order 11 on $\Pic(X)\otimes \bbQ$, we get (iv) from (iii).
\end{proof}

\begin{corollary}\label{jac} Let $X$ be a K3 surface over an algebraically closed field $k$
of characteristic $11$. Assume that $X$ admits an automorphism $g$
of order $11$. Then $X$ is isomorphic to a torsor of one of the
elliptic surfaces $X_\varepsilon$. The order of this torsor in the
Shafarevich-Tate group of torsors is equal to $1$ or $11$.
\end{corollary}

\begin{proof} Let $f_J:J\to \bbP^1$ be the Jacobian fibration of the
elliptic fibration $f:X\to \bbP^1$ defined by the $g$-invariant
elliptic pencil. Let $J^o$ be the open subset of $J$ whose
complement is the set of singular fibres of $f_J$. We know that
the fibres of $f$ are irreducible. By a result of M. Raynaud, this
allows us to identify  $J^o$ with the component
$\textbf{Pic}_{X/\bbP^1}^0$ of the relative Picard scheme of
invertible sheaves of degree 0 (see \cite{CD}, Proposition 5.2.2).
The automorphism $g$ acts naturally on the Picard functor and
hence on $J^o$. Since $J$ is minimal, it acts biregularly on $J$.
This action preserves the elliptic fibration on $J$ and defines an
automorphism of order $11$ on the base. This implies that there
exists an $C_{11}$-equivariant isomorphism of elliptic surfaces
$J$ and $X_\varepsilon$.

The assertion about the order of the torsor follows from the
existence of a section or an 11-section of $f$. In fact, let $Y$
be a nonsingular relatively minimal model of the elliptic surface
$X/(g)$ with the elliptic fibration induced by $f$. It is a
rational elliptic surface. Let $F_0$ be the $g$-invariant fibre of
$f$ over a point $s_0\in \bbP^1$. The singular fibres of the
elliptic fibration $f':Y\to \bbP^1$ over $\bbP^1\setminus\{s_0\}$
are either two irreducible nodal fibres ($\varepsilon\neq 0$) or
one cuspidal irreducible fibre ($\varepsilon = 0$). The standard
argument in the theory of elliptic surfaces shows that the fibre
of $f'$ over $s_0$ is either of type $\tilde{E}_8$ or
$\tilde{D}_8$. This fibre is not multiple if and only if $f'$ has
a section. The pre-image of this section is a section of $f$
making $X$ the trivial torsor. A singular fibre of additive type
can be multiple only if the characteristic is positive, and the
multiplicity $m$ must be equal to the characteristic (see
\cite{CD}, Proposition 5.1.5). In this case an exceptional curve
of the first kind on $Y$ is a $m$-section. The pre-image of this
multi-section on $X$ is a $m$-section, in our case an
$11$-section.
\end{proof}

\begin{remark} Note that, even in the case $X = X_\varepsilon$, the $g$-invariant fibration
may be different from the standard elliptic fibration. In other
words a non-trivial torsor of an elliptic surface could be
isomorphic to the same surface. This strange phenomenon could
happen only in positive characteristic and only for torsors of
order divisible by the characteristic. We do not know an example
where this strange phenomenon really occurs. In Kond\=o's example,
the $g$-invariant elliptic fibration for an element $g$ of order
11 in $G=M_{11}$ or $M_{22}$ may have a section (but no
$g$-invariant section!). If this happens, it is isomorphic to the
standard elliptic fibration and hence $g$ is conjugate to $\tau
g_\varepsilon$ as we have seen in Proposition \ref{surface}.
\end{remark}

\begin{lemma}\label{lemma} Suppose $p = 11$. Then there is a finite subgroup
$K_\varepsilon$ of $\Aut(\xvar)$ satisfying the following property:
\begin{itemize}
\item[(i)] $K_\varepsilon$ leaves invariant both the standard elliptic fibration
of $\xvar$ and the zero section $S_\varepsilon$ which is the
section at infinity.
\item[(ii)] $K_0\cong \GU_2(11)/(\pm I)\cong L_2(11):12$ and  $K_1\cong 11:4$, where the
first factor in the semi-direct product is a symplectic subgroup and
the second factor a non-symplectic subgroup.
\item[(iii)] The image of $K_\varepsilon$ in  $\Aut(\bbP^1)$ is equal to
the subgroup $\Aut(\bbP^1, V(\Delta_\varepsilon))$ which leaves
the set $V(\Delta_\varepsilon)$ invariant.
\item[(iv)]  $\Aut(\bbP^1, V(\Delta_0)) \cong
  \PGU_2(11) \cong L_2(11).2$ and  $\Aut(\bbP^1,
 V(\Delta_\varepsilon)) \cong 11:2$ if $\varepsilon\neq 0$.
\end{itemize}
\end{lemma}

\begin{proof} Assume $\varepsilon = 0$. After a linear change of variables
$$t_0 = \alpha^{11}t'_0+\alpha t'_1,\quad t_1 = t'_0+t'_1,$$
where $\alpha\in \bbF_{{11}^2}\setminus \bbF_{11}\subset k^*$, we can
transform the polynomial $t_0t_1^{11}-t_0^{11}t_1$ to the form
$\lambda t_0^{12}+\mu t_1^{12}$. After scaling, it becomes of the
form $f = t_0^{12}+ t_1^{12}$. Now notice that this equation
represents a hermitian form over the field $\bbF_{{11}^2}$, hence the
finite unitary group $\GU_2(11)$ leaves the polynomial $f$
invariant. The group $\GU_2(11)$ acts on the surface
\begin{equation}\label{neweq}
X_0 \cong V(y^2+x^3+t_0^{12}+t_1^{12})
\end{equation}
in an obvious way, by acting on the variables
$t_0,t_1$ and identically on the variables $x,y$. Note that
$$(t_0, t_1, x, y)=(\lambda t_0, \lambda t_1, \lambda^4 x,
\lambda^6 y)$$
in $\bbP(1,1,4,6)$ for all $\lambda \in k^*$. In particular
$(t_0, t_1, x, y)=(-t_0, -t_1, x, y)$, so $-I\in
\GU_2(11)$ acts trivially on $X_0$. Note also that $\SU_2(11)$ and hence $\PSU_2(11)$
acts symplectically on $X_0$. The action of $\PSU_2(11)$ is faithful
because it is a simple group.  Take $K_0=\GU_2(11)/(\pm I)$ and
consider the homomorphism $$\det: K_0\to (\bbF_{{11}^2})^*.$$
It is known that
 $$U_2(11) = \PSU_2(11)\cong \PSL_2(\bbF_{11})=L_2(11).$$
If $A\in \GU_2(11)$, then
$(\det A)^{12}=(\det A)(\overline{\det A})=\det A^t\overline{A}
=\det I=1$, so the image of det is a cyclic group of order dividing
12. On the other hand, if $\zeta\in
\bbF_{{11}^2}$ is a 12-th root of unity, the unitary matrix
\[\begin{pmatrix}1&0\\
0&\zeta\end{pmatrix}\]
generates an order 12 subgroup of $K_0$, which acts on $X_0$ non-symplectically.
This proves (i) and (ii).

We know that the group $\GU_2(11)$ leaves the polynomial $f$
invariant. Thus its image $\PGU_2(11)$ in $\Aut(\bbP^1)$ must coinside
with  $\Aut(\bbP^1, V(\Delta_0))$. It is known
 that $\PGU_2(11)$ is a
maximal subgroup in the permutation group $\frakS_{12}$ and
$\PGU_2(11) \cong \PGL_2(\bbF_{11})\cong L_2(11).2$, a non-split
extension. The quotient group
is generated by the image of the automorphism $:(t_0, t_1)\to (t_0,
\zeta t_1)$, where  $\zeta\in\bbF_{{11}^2}$ is a 12-th root of unity.
This proves (iii) and (iv).

Assume $\varepsilon \neq 0$. An element of $\PGL_2(k)$ leaving
$V(\Delta_\varepsilon)$ invariant must either leave all factors of
$\Delta_\varepsilon$ from \eqref{for2} invariant or interchange
the second and the third factors. It can be seen by computation
that the group $\Aut(\bbP^1,
 V(\Delta_\varepsilon))$ is generated by the following 2 automorphisms
$$e(t_0, t_1)=(t_0, t_1+t_0), \quad i(t_0, t_1)=(t_0, -t_1+bt_0)$$
where $b$ is a root of $b^{11}-b+3\varepsilon^3=0$. The order of
$e$ (resp. $i$) is 11 (resp. 2) and $i$ normalizes $e$. We see
that they lift to automorphisms of $X_\varepsilon$
$$\tilde{e}(t_0, t_1, x, y)=(t_0, t_1+t_0, x, y), \,
\tilde{i}(t_0, t_1, x, y)=(t_0, -t_1+bt_0, -x+3\varepsilon t_0^4,
\sqrt{-1}y)$$ and we take $K_\varepsilon=(\tilde{e}, \tilde{i})$.
Clearly $\tilde{i}$ is non-symplectic of order 4 and normalizes
$\tilde{e}$ which is symplectic of order 11, and both leave
invariant the zero section $S_\varepsilon$.
\end{proof}

\begin{remark}\label{goto} The equation \eqref{neweq} makes $X_0$  a
  weighted Delsarte surface according to the definition in
  \cite{Goto}. It follows from loc.cit. that $X_0$ is a supersingular
  surface with Artin invariant $\sigma = 1$. It follows from the
  uniqueness of such surface  that $X_0$ is also isomorphic to the Fermat quartic
$$x_0^4+x_1^4+x_2^4+x_3^4 = 0,$$
 the Kummer surface associated to the product of supersingular elliptic curves,
and the modular elliptic surface of level 4 (see \cite{Shioda}).
We do not know whether the surface $X_\varepsilon$,
$\varepsilon\neq 0$, is supersingular. By Lemma \ref{lem1}, we
know that $\rank~\Pic(X_\varepsilon)=2$, 12 or 22.
\end{remark}

\begin{definition}\label{H} The subgroup
  $K_\varepsilon\subset\Aut(\xvar)$ from Lemma \ref{lemma}
contains a symplectic subgroup leaving invariant the standard
elliptic fibration of $\xvar$, isomorphic to $L_2(11)$ if
$\varepsilon = 0$ and  to $C_{11}$ if $\varepsilon = 1$. Denote
this subgroup by $H_{\varepsilon}$. It leaves invariant the zero
section $S_\varepsilon$ of the elliptic fibration.
\end{definition}

The group $H_{\varepsilon}$ acts on the base curve $\bbP^1$ and we
have a homomorphism
$$\pi: H_{\varepsilon}\to \Aut(\bbP^1, V(\Delta_\varepsilon)),$$
which is an embedding. The image $\pi(H_{\varepsilon})$ is equal to
the unique index 2 subgroup of
$\Aut(\bbP^1, V(\Delta_\varepsilon))$.

\begin{proposition}\label{Hmax} Let $G$
be a finite group of symplectic automorphisms of the surface
$\xvar$ leaving invariant the standard elliptic fibration of
$\xvar$. Let
$$\psi: G\to \Aut(\bbP^1, V(\Delta_\varepsilon))$$
be the natural homomorphism. Then $\psi$ is an embedding. If in
addition $G$ is wild and leaves invariant the zero section
$S_\varepsilon$, then $G$ is contained in $H_{\varepsilon}$.
\end{proposition}

\begin{proof}
Let $\alpha\in \Ker(\psi)$. Then $\alpha$ acts trivially on the
base curve. Since $p>3$,  $\alpha$ being symplectic must be a
translation by a torsion section. It is known that there is no
$p$-torsion in the Mordell-Weil group of an elliptic K3 surface if
the characteristic $p>7$ (\cite{DK2}), and there are no other
torsion sections because no symplectic automorphism of order
coprime to $p$ can have more than 8 fixed points (Theorem 3.3
\cite{DK2}), while the fibration has 12 or 23 singular fibres.
Hence $\alpha$ must be the identity automorphism. This proves that
$\psi$ is an embedding.

If $\psi$ is surjective, then $\#G=2.\#L_2(11)$ or $2.11$, which
cannot be the order of a wild symplectic group in characteristic
11, by Proposition \ref{mathieu} and Lemma \ref{orders}. Thus
$\psi$ is not surjective. From this we see that if $G$ is wild,
then $\psi(G)$ is contained in the unique index 2 subgroup
$\pi(H_{\varepsilon})$ of $\Aut(\bbP^1, V(\Delta_\varepsilon))$.
If an element $\alpha\in G$ and an element $h\in H_\varepsilon$
have the same image in $\Aut(\bbP^1, V(\Delta_\varepsilon))$, then
$\alpha h^{-1}$ is a translation by a section. If $\alpha$ leaves
invariant the zero section $S_\varepsilon$, so does $\alpha
h^{-1}$, hence $\alpha h^{-1}$ is the identity automorphism. This
proves the second assertion.
\end{proof}

\section{A Mathieu representation}
From now on $X$ is a K3 surface over an algebraically closed field of
characteristic $p = 11$ and $G$ a group of symplectic automorphisms
of $X$ of order divisible by $11$.

\begin{lemma}\label{lemm2} Let $S$ be a normal projective rational
  surface with an isolated singularity $s$. Then
$$e_c(S\setminus \{s\}) \ge 2,$$
where $e_c$ denotes the l-adic Euler-Poincar\'e characteristic with compact support.
\end{lemma}

\begin{proof} Let $f:S'\to S$ be a minimal resolution of $S$. Let $E$
  be the reduced exceptional divisor. Then $e_c(E) = 1-b_1(E)+b_2(E)
  \le 1+b_2(E)$. Since the intersection
matrix of irreducible components of $E$ is negative definite, we have
$b_2(S') \ge 1+b_2(E)$. This gives
$$e_c(S\setminus \{s\}) = e_c(S'\setminus E) = e_c(S')-e_c(E) \ge 2+b_2(S')-(1+b_2(E)) \ge 2.$$
\end{proof}

\begin{lemma}\label{lemm3} Let $g$ be an automorphism of $X$ of
order $11$. Assume that $X^g$ is a point. Then the cyclic group $(g)$
is not contained in a larger symplectic cyclic subgroup of $\Aut(X)$.
\end{lemma}

\begin{proof} Let $H=(h)$ be a symplectic cyclic subgroup of $\Aut(X)$ containing
  $(g)$. Write $\#H=11r$ and $g=h^r$. Without loss of generality,
we may assume that $r$ is a prime, and by Theorem 3.3  of \cite{DK2},
may further assume that $r=2,3,5,7$, or $11$.

Assume $r\neq 11$. Let $f=h^{11}$. Then $f$ is symplectic of order
$r=2,3,5$, or $7$. By Theorem 3.3  of \cite{DK2}, $X^f$ is a
finite set of points of cardinality $< 11$. The order 11
automorphism $g$ acts on $X^f$, hence acts trivially. Thus
$X^f\subset X^g$, but $\#X^f\ge 3$, a contradiction. Thus $r=11$.

Let $x\in X$ be the fixed point of $g$, and $y\in X/(g)$ be its
image. Let $V = X/(g)\setminus \{y\}$. We claim that the quotient
group $\bar{H}=H/(g)$ acts freely on $V$. To see this, suppose
that $h(z)=g^i(z)$ for some point $z\in X$, some $g^i\in (g)$.
Then $g(z)=h^{11}(z)=g^{11i}(z)=1_X(z)=z$, so $z=x$. This proves
the claim.

By Lemma \ref{lem1}, for any $l\ne 11$, $\dim H^2_{\rm
et}(X,{\bbQ}_l)^g=2$. This implies that $$\Tr (g^*|H^2_{\rm
et}(X,{\bbQ}_l))=0.$$ By the Trace formula of S. Saito
\cite{Saito},
$$l_x(g)=\Tr (g^*|H^*_{\rm et}(X,{\bbQ}_l))=2,$$
where $l_x(g)$ is the intersection index of the graph of $g$ with
the diagonal at the point $(x,x)$. The formula of Saito
(\cite{Saito}, Theorem 7.4, or \cite{DK2}, Lemma 2.8) gives
$e_c(V)= 3$. Since the group $H/(g)$ acts freely on $V$,
$e_c(V/\bar{H}) = 3/\#\bar{H}.$ Applying Lemma \ref{lemm2} to the
surface $S = X/H$, we obtain that $\bar{H}$ is trivial.
\end{proof}

\begin{lemma} \label{order} Let $G$ be a finite group of symplectic
automorphisms of a K3 surface $X$ over an algebraically closed
field of characteristic $p=11$. Then $${\rm ord}(g) \in
\{1,2,3,4,5,6,7,8,11\}$$ for all $g\in G$.
\end{lemma}

\begin{proof} If the order ${\rm ord}(g)$ of $g\in G$ is coprime to the
characteristic $p=11$, then by Theorem 3.3  of \cite{DK2} $${\rm
ord}(g)\in \{1,\ldots,8\}.$$ It remains to show that $G$ cannot
contain any element of order $11r, r>1$. Assume the contrary, and
let $h\in G$ be an element of order $11r$. We may assume that $r$
is a prime and hence $r=2,3,5,7$, or $11$. Let $g=h^r$ and $f =
h^{11}$. We see that $g$ is of order 11. By Lemma \ref{lemm3},
$X^g$ cannot be a point, hence must be a cuspidal curve
(\cite{DK1}). Denote this curve by $F$. It is easy to see that $F$
is $h$-invariant, i.e. $h(F)=F$.

Assume $r=11$. Then $h$ acts on the base curve $\bbP^1$ of the
pencil $|F|$ faithfully, however, using the Jordan canonical form
we see that $\bbP^1$ does not admit an automorphism of order
$11^2$.

Next, assume that $r=2,3,5,7$. By Theorem 3.3  of \cite{DK2},
$$3\le \#X^{f}\le 8.$$
Since $r$ is prime to $11$,
$$X^h=X^{f}\cap X^{g}.$$
Clearly $g$ acts on the finite set $X^{f}$, and this action cannot
be of order $11$. This means that $g$ acts trivially on $X^{f}$,
i.e. $X^{f}\subset X^{g}=F$. Thus
$$X^h= X^{f}.$$
This means that $h$ acts on $F$ with $\#X^f$ fixed points. But no
nontrivial action on a rational
 curve can fix more than 2 points. A contradiction.
\end{proof}

\medskip
A Mathieu representation of a finite group $G$ is a 24-dimensional
representation  on a vector space $V$ over a field of characteristic
zero with character
$$\chi(g)=\epsilon({\rm ord}(g)),$$ where
\begin{equation}\label{mumu}
\epsilon(n)=24(n\prod_{p|n}(1+{1\over p}))^{-1}, \quad \epsilon(1)=24.
\end{equation}
The number
\begin{equation}\label{mu}
\mu(G) = \frac{1}{\# G}\sum_{g\in G}\epsilon(\ord(g))
\end{equation}
is equal to the dimension of the subspace $V^G$ of $V$. Here $V^G$ is
the linear subspace  of vectors fixed by $G$.
The natural action of a finite group $G$ of symplectic
automorphisms of a complex K3 surface on the singular cohomology
$$H^*(X,{\bbQ})=\oplus_{i=0}^4H^i(X,{\bbQ})\cong {\bbQ}^{24}$$ is a Mathieu representation with
$$\mu(G) = \dim H^*(X,{\bbQ})^G \ge 5.$$ From this Mukai deduces that $G$ is
isomorphic to a subgroup of $M_{23}$ with at least 5 orbits.
In positive characteristic, if $G$ is
wild, then the
formula for the number of fixed points is no longer true and the
representation of $G$ on the $l$-adic cohomology, $l\ne p$,
$$H^*_{\rm et}(X,{\bbQ}_l)=\bigoplus_{i=0}^4H^i_{\rm
  et}(X,{\bbQ}_l)\cong {\bbQ}_l^{24}$$ is not Mathieu in general.
However, in our case we have the following:

\begin{proposition} \label{mathieu} Let $G $ be a finite group
  acting symplectically on a K3 surface $X$ over a field of
  characteristic $11$. Then the representation of $G$ on $H^*_{\rm
  et}(X,\bbQ_l)$, $l\ne 11$,
is a Mathieu representation with $\dim H^*_{\rm et}(X,\bbQ_l)^G \ge 3$.
\end{proposition}

\begin{proof} Note that $\rank~\Pic(X)^G\ge 1$, and the second
assertion follows. It remains to prove that the representation is
Mathieu.  By Lemma \ref{order}, it is enough to show this for
automorphisms of order $11$. Let $g\in G$ be an element of order
$11$. We need to show that the character $\chi(g)$ of the
representation on the $l$-adic cohomology $H^*_{\rm
et}(X,{\bbQ}_l)$ is equal to $\epsilon(11)=2$. Since
$$\chi(g)=\Tr (g^*|H^*_{\rm et}(X,{\bbQ}_l)),$$
it suffices to show that $\Tr (g^*|H^2_{\rm et}(X,{\bbQ}_l))=0,$
or $\dim H^2_{\rm et}(X,{\bbQ}_l)^g=2$.  Now the result follows
from Lemma \ref{lem1}.
\end{proof}

\section{Determination of Groups}
In this section we determine all possible finite groups which may
act symplectically and wildly on a K3 surface in characteristic
11. We use only purely group theoretic arguments.

\begin{proposition}\label{bd} Let $G$ be a finite group having
a Mathieu representation over $\bbQ$ or over $\bbQ_l$ for all prime $l\ne 11$. Then
$$\#G=2^a.3^b.5^c.7^d.11^e.23^f$$
for $a\le 7, b\le 2, c\le 1, d\le 1,
e\le 1, f\le 1$.
\end{proposition}

\begin{proof}
If the representation is over $\bbQ$, this is the theorem of Mukai
(\cite{Mukai} (Theorem (3.22)).
In his proof, Mukai uses at several places the fact that
the representation is over $\bbQ$.  The only essential place where
he uses that the representation is over $\bbQ$ is Proposition
(3.21), where $G$ is assumed to be a 2-group containing a maximal
normal abelian subgroup $A$ and the case of $A =(\bbZ/4)^2$ with
$\#(G/A)\ge 2^4$ is excluded by using that a certain 2-dimensional representation of
the quaternion group $Q_8$ cannot be defined over $\bbQ$. We use that $G$
also admits a Mathieu representation on 2-adic cohomology, and it
is easy to see that the representation of  $Q_8$ cannot be defined
over $\bbQ_2$.
\end{proof}

The following lemma is of purely group theoretic nature and its proof
follows an argument employed by S. Mukai \cite{Mukai}.

\begin{lemma}\label{orders} Let $G$ be a finite group having
a Mathieu representation over $\bbQ$ or over $\bbQ_l$ for all prime $l\ne 11$.
Assume $\mu(G)\ge 3$. Assume that $G$ contains
an element of order $11$, but no elements of order $>11$.
Then the order of $G$ is equal to one of the following:
$$11,\quad 5.11, \quad 2^2.3.5.11, \quad 2^4.3^2.5.11,
\quad 2^7.3^2.5.7.11.$$
\end{lemma}

\begin{proof}
Since $G$ has no elements of order 23, by  Proposition \ref{bd}, we have
$$\#G=2^a.3^b.5^c.7^d.11, \, (a\le 7, b\le 2, c\le 1, d\le 1).$$
Let $S_q$ be a $q$-Sylow subgroup of $G$ for $q=5,7$ or $11$. Then
$S_q$ is cyclic and its centralizer coincides with $S_q$. Let $N_q$ be
the normalizer of $S_q$. Since $G$ does not contain elements of order
$5k, 7k, 11k$ with $k > 1$,  the index $m_q:=[N_q:S_q]$ is a divisor
of $q-1$. Since it is known that the dihedral groups $D_{14}$ and $D_{22}$ do not admit a
Mathieu representation, we have $m_7=1$ or 3, and $m_{11}=1$ or 5.
Let $a_n$ be the number of elements of order $n$ in $G$. Then
$a_q=\frac{\#G(q-1)}{qm_q}.$
As in \cite{Mukai}, we have
\begin{equation}\label{muG}
\mu(G)=\frac{1}{\#G}\sum\epsilon(n)a_n
=8+\frac{1}{\#G}(16-2a_3-4a_4-4a_5-6a_6-5a_7-6a_8-6a_{11}).
\end{equation}

\medskip\noindent
Case 1. The order of $G$ is divisible by 7.

\medskip
The formula \eqref{muG} gives
\begin{equation}\label{muG2}
\mu(G)\le 8+\frac{16}{\#G}-\frac{30}{7m_7}-\frac{60}{11m_{11}}.
\end{equation}
Since $\mu(G) \ge 3$, the numbers $m_{11}$ and $m_{7}$ must be greater than 1.

Assume $m_{11}=5, m_{7}=3$. Then $\#G$ is divisible by 5, and the formula \eqref{muG} gives
\begin{equation}\label{muG3}
\mu(G)\le 8+\frac{16}{\#G}-\frac{16}{5m_5}-\frac{10}{7}-\frac{12}{11}.
\end{equation}
If $m_5=1$, then this inequality gives $\mu(G)<3$. If $m_5=2$, then
the number of $q$-Sylow subgroups is equal to $2^{a-1}.3^b.7.11$,
$2^{a}.3^{b-1}.5.11$, $2^{a}.3^b.5.7$ for $q=5,7,11$ respectively. Taking $q = 5$ and applying
Sylow's theorem, we get  $a-b\equiv 0 \mod 4$. Since $1\le a\le 7$,
$1\le b\le 2$, the only solutions are $(a,b) =  (5,1),
(6,2)$. However, neither $2^5.5.11$ nor
$2^6.3.5.11$ is congruent to 1 modulo 7.

If $m_5=4$, then
the number of $q$-Sylow subgroups is equal to $2^{a-2}.3^b.7.11$,
$2^{a}.3^{b-1}.5.11$, $2^{a}.3^b.5.7$ for $q=5,7,11$ respectively. A similar argument, shows that  $a-b \equiv 1 \mod 4$ and the possible order is $2^7.3^2.5.7.11$.

\medskip\noindent
Case 2. The order of $G$ is divisible by 5, but not by 7.

\medskip
The formula \eqref{muG} gives
\begin{equation}\label{muG4}
\mu(G)\le 8+\frac{16}{\#G}-\frac{16}{5m_5}-\frac{60}{11m_{11}}.
\end{equation}
Assume that $m_{11}=1$. Then this
inequality gives $\mu(G)<3$.

Assume $m_{11}=5$.
If $m_5=1$, then the number of $q$-Sylow subgroups is equal to $2^{a}.3^b.11$,
$2^{a}.3^b$ for $q=5, 11$ respectively. By
Sylow's theorem,  $a-b\equiv 0 \mod 4$,
$a+8b \equiv 0 \mod 10$. This system of
congruences has only one solution $a=b=0$ in the range $a\le 7$, $b\le
2$. This gives the possible order $5.11$.

If $m_5=2$, then
the number of $q$-Sylow subgroups is equal to $2^{a-1}.3^b.11$,
$2^{a}.3^b$ for $q=5, 11$ respectively. By
Sylow's theorem,  $a-b\equiv 1 \mod 4$,
$a+8b \equiv 0 \mod 10$. This system
has only one solution $a=2, b=1$ in the range $1\le a\le
7$, $b\le 2$. This gives the possible order $2^2.3.5.11$.

If $m_5=4$, then
the number of $q$-Sylow subgroups is equal to $2^{a-2}.3^b.11$,
$2^{a}.3^b$ for $q=5, 11$ respectively. By
Sylow's theorem,  $a-b\equiv 2 \mod 4$,
$a+8b \equiv 0 \mod 10$. This system has only one solution
$a=4, b=2$ in the range $2\le a\le 7$,
$b\le 2$. This gives the possible order $2^4.3^2.5.11$.

\medskip\noindent
Case 3. The order of $G$ is divisible by neither 5 nor 7.

\medskip
In this case $m_{11}\ne 5$, and hence $m_{11}=1$. Thus the formula \eqref{muG} gives
\begin{equation}\label{muG5}
\mu(G)\le 8+\frac{16}{\#G}-\frac{60}{11}.
\end{equation}
The number of $11$-Sylow subgroups is equal to $2^{a}.3^b$. By
Sylow's theorem, $a+8b \equiv 0 \mod 10$. This
congruence has 3 solutions $(a,b)=(0,0)$, $(2, 1)$, $(4, 2)$ in the range $a\le 7$, $b\le
2$. The first gives the possible order $11$.
In the second and the third case, the inequality \eqref{muG5}
gives $\mu(G)<3$.
\end{proof}

\begin{proposition}\label{simple} In the situation of the previous
lemma, $G$ is isomorphic to one of the following groups:
$$C_{11},\quad 11:5, \quad L_2(11), \quad M_{11}, \quad M_{22}.$$
\end{proposition}

\begin{proof} By Lemma \ref{orders}, there are five possible orders for
$G$
\begin{equation}\label{five}
11,\quad 5.11, \quad 2^2.3.5.11, \quad 2^4.3^2.5.11, \quad
2^7.3^2.5.7.11.
\end{equation}
In the first two cases, the assertion is obvious.
The remaining possible 3 orders are exactly the orders of the 3 simple
groups given in the assertion. The theory of finite simple groups shows that
there is only one simple group of the order in each of these cases.

Assume the last 3 cases. It suffices to show that $G$ is simple.

Let $K$ be a proper normal subgroup of $G$ such that $G/K$ is simple. If $\#K$ is not divisible
by 11, then an order 11 element of $G$ acts on the set $\Syl_q(K)$ of
$q$-Sylow subgroups of $K$. Since
$\#\Syl_q(K)$ is not divisible by 11 for any prime $q$ dividing $\#K$,
the order 11 element $g$ must normalize a $q$-Sylow subgroup of $K$.
If one of the numbers $q = 3,5,$ or $7$ divides $\#K$, then $g$
centralizes an element of one of these orders. This
contradicts  the assumption that $G$ does not contain an element of
order $>11$. If $q =2$ divides $\#K$, then a
$2$-Sylow subgroup of $K$ is of order $\le 2^7$, and hence $g$
centralizes an element of order 2, again a
contradiction.
So, we may assume that
$11|\#K$. If $\#K=11$, then an order 2 element of $G$ normalizes
$K$. Neither a cyclic group of order 22 nor a dihedral group of order
22 has a Mathieu representation,  so
$\#K>11$. If $K\cong 11:5$, then an order 2 element of $G$ normalizes
the unique $11$-Sylow subgroup of $K$, again a contradiction. If
$\#K$ is one of the remaining three possibilities, then  the group $G/K$
is of order $2^5.3.7$ or $2^3.7$ or $2^2.3$. In the first case an
order 7 element of $G$ normalizes, hence centralizes a Sylow
11-subgroup of $K$, again a contradiction. Obviously in the
other two cases $G/K$ cannot be simple. This proves that $G$ is simple.
\end{proof}

\begin{corollary}\label{normalizer} Let $G $ be a finite group
  acting symplectically and wildly on a K3 surface $X$ over a field of
  characteristic $11$. Let $g$ be an element of order $11$ in $G$.
  Then the normalizer of $(g)$ in $G$ must be isomorphic to $11:5$ if $\#G>11$.
\end{corollary}

\section{Proof of the Main Theorem}

In this section we complete the proof of Theorem \ref{main}
announced in Introduction. It remains to prove the assertion (ii).

\begin{lemma}\label{varneq0} Assume $\varepsilon\neq 0$. Let $G\subset \Aut(X_\varepsilon)$
be a finite wild symplectic subgroup. If an element $g\in G$ of
order $11$ leaves invariant the standard elliptic fibration with a
$g$-invariant section, then $G=(g)\cong C_{11}$ and $G$ is
conjugate to $H_\varepsilon=(g_\varepsilon)$. In particular,
$H_\varepsilon$ is a maximal finite wild symplectic subgroup of
$\Aut(X_\varepsilon)$.
\end{lemma}

\begin{proof} Since $g$ leaves a section invariant, it must be a
conjugate to $g_\varepsilon$. So up to conjugation, we may assume
that $g$ leaves the zero section $S_\varepsilon$ invariant. Thus
$g=g_\varepsilon$ by Proposition \ref{Hmax}.

Suppose $G>(g)$. Let $N$ be the normalizer of $(g)$ in $G$. Then
$N \cong 11:5$ by Corollary \ref{normalizer}.

Claim that $N$ leaves invariant the standard elliptic pencil
$|F|$. It is enough to show that $h(F_0)=F_0$ for any $h\in N$,
where $F_0=X^g$ is a cuspidal curve in $|F|$. In fact, for any
$x\in F_0$, we have $h(x)=hg(x) = g^ih(x)$ for some $i$, so
$h(x)\in X^{(g)}=F_0$, which proves the claim.

Next, claim that $N$ leaves invariant the zero section
$S_\varepsilon$. In fact,
$h(S_\varepsilon)=hg(S_\varepsilon)=g^ih(S_\varepsilon)$, so $(g)$
leaves invariant $h(S_\varepsilon)$, and hence
$h(S_\varepsilon)=S_\varepsilon$ as $g$ cannot leave invariant two
distinct sections by Lemma \ref{lem1} (iii).

Now Proposition \ref{Hmax} gives a contradiction. Hence, $G=(g)$.
\end{proof}

\begin{lemma}\label{var0} Let $G\subset \Aut(X_0)$
be a finite wild symplectic subgroup, isomorphic to $L_2(11)$. If
an element $g\in G$ of order $11$ leaves invariant both the
standard elliptic fibration and a section, then $G$ is conjugate
to $H_0$. In particular, if $G$ contains $g_0$ then $G=H_0$.
\end{lemma}

\begin{proof} Replacing $G$ by a conjugate subgroup in $\Aut(X_0)$, we may assume
that $g$ leaves invariant both the standard elliptic fibration and
the zero section $S_0$, i.e. $g=g_0$. We need to prove that
$G=H_0$.

Let $|F|$ be the standard elliptic fibration. Then $g(S_0)=S_0$
and $X^g=F_0$, a cuspidal curve in $|F|$.

Let $N$ be the normalizer of $(g)$ in $G$. Then $N \cong 11:5$.
The same argument as in the proof of Lemma \ref{varneq0} shows
that $N$ leaves invariant both the standard elliptic pencil $|F|$
and the zero section $S_0$. By Proposition \ref{Hmax}, $N\subset
H_0$.

We have $N\subset G\cap H_0$. Suppose $G\cap H_0=N$. Consider the
$G$-orbit of the divisor class $[F]\in\Pic(X_0)$,
$$G([F])=\{ h([F])\in\Pic(X_0)| h\in G\}.$$
Clearly $N$ acts on it. Note $$\#G([F])=[G:N]=12.$$ Thus $G([F])$
is the set of 12 different elliptic fibrations with a section. The
automorphism $g$ cannot leave invariant an elliptic fibration
other than $|F|$, hence fixes $[F]$ and has 1 orbit on the
remaining 11 elliptic fibrations, which we denote by
$[F_1],\cdots, [F_{11}]$.

Recall that $H_0$ leaves invariant the zero section $S_0$. The
three divisor classes
$$[F], \quad \sum_{j=1}^{11}[F_j], \quad [S_0]$$
are $N$-invariant, and their intersection matrix is given as follows:
\[\begin{pmatrix}0&11m&1\\
11m&110m&11b\\
1&11b&-2\end{pmatrix}\] where $m = F\cdot F_i, \ b = S_0\cdot F_i,
i \ge 1$. Its determinant is equal to $$242(m^2+bm)-110m,$$ which
cannot be 0 for any positive integers $m$ and $b$. This implies
that
$$\mu(N) = 2+\rank~\Pic(X_0)^{N}\ge 5,$$
a contradiction to the equality $\mu(N)=4$. This proves that $N$
is a proper subgroup of $G\cap H_0$. Since $N$ is a maximal
subgroup of $G$, we have $G=H_0$.
\end{proof}

Note that $\mu(M_{11})=\mu(M_{22})=3$ and $\mu(L_2(11))=4$. Note also
that $L_2(11)$ is isomorphic to a maximal subgroup of both $M_{11}$
and $M_{22}$.

The following proposition completes the proof of Theorem
\ref{main} (ii).

\begin{proposition}\label{fin} Let $G\subset \Aut(X_0)$
be a finite wild symplectic subgroup. Assume that $G\cong M_{11}$
or $M_{22}$. Then no conjugate of $G$ in $\Aut(X_0)$ contains the
automorphism $g_0$ given by \eqref{form2}. In other words, no
element of $G$ of order $11$ can leave invariant both the standard
elliptic fibration and a section. In particular, $H_0$ is a
maximal finite wild symplectic subgroup of $\Aut(X_0)$.
\end{proposition}

\begin{proof} Suppose that a conjugate of $G$ contains
$g_0$. Replacing $G$ by the conjugate, we may assume that $g_0\in
G$.

Let $K$ be a subgroup of $G$ such that $g_0\in K\subset G$ and
$K\cong L_2(11)$. Then by Lemma \ref{var0}, $K=H_0$. Thus $$g_0\in
H_0\subset G.$$ Since $H_0\cong L_2(11)$ is a maximal subgroup of
$G$, its normalizer subgroup $N_G(H_0)$ coincides with $H_0$.

Let $|F|$ be the standard elliptic fibration on $X_0$, and $S_0$
the zero section. Then $g(S_0)=S_0$ and $X^g=F_0$, a cuspidal
curve in $|F|$. Furthermore, both the section $S_0$ and the
elliptic pencil $|F|$ are $H_0$-invariant (see Definition
\ref{H}).

Consider the $G$-orbit of the divisor class $[F]$,
$$G([F])=\{ h([F])\in\Pic(X_0)| h\in G\}.$$
Consider the action of $H_0$ on it. By Proposition \ref{Hmax}, the
stabilizer subgroup $G_{[F]}$ of $[F]$ coincides with $H_0$. The
automorphism $g_0$ cannot leave invariant two different elliptic
fibrations, hence fixes $[F]$ and has orbits on the set
$G([F])\setminus \{[F]\}$ of cardinality divisible by 11. This
implies that $H_0$ fixes $[F]$ and has orbits on the set
$G([F])\setminus \{[F]\}$ of cardinality divisible by 11. Write
$$G([F])=\{[F]=[F_0],\, [F_1],\, [F_2],\, ...,\, [F_{r-1}]\}$$
where $r=\#G([F])=[G:H_0]$. Let $$\calO_1\cup \calO_2\cup
...\cup\calO_s$$ be the orbit decomposition of the index set
$\{1,\, 2,\, ...,\, r-1\}$ by the action of $H_0$. Since $H_0$
fixes $[F]$ and acts transitively on each $\calO_i$, the
intersection number $F\cdot F_t$ is constant on the orbit
$\calO_i$ containing $t$, i.e. $F\cdot F_t=m_i$ for all $t\in
\calO_i$. Note that the divisor $$\calF=\sum_{j=0}^{r-1}F_j$$ is
$G$-invariant, and
\begin{equation}\label{comp1}
\calF^2=(\sum_{j=0}^{r-1}F_j)^2=rF_0\cdot\sum_{j=0}^{r-1}F_j=r\sum_{i=1}^sm_i\#\calO_i.
\end{equation}

Next recall that $H_0$ leaves invariant the zero section $S_0$.
Similarly we consider the $G$-orbit of the divisor class $[S_0]$
$$G([S_0])=\{ h([S_0])\in\Pic(X_0)| h\in G\}.$$
Let $G_0$ be the stabilizer subgroup of $[S_0]$. Since it contains
$H_0$ and $H_0$ is maximal in $G$, we obtain that $G_0 = H_0$ or
$G_0 = G$.

Assume $G_0 = H_0$. Then all stabilizers are conjugate to $H_0$.
Similarly as above we claim that $g_0\in H_0$ fixes no elements of
$G([S_0])$ other than $[S_0]$. If $g_0h(S_0)=h(S_0)$ for some
$h\in G$, Then $g_0\in hH_0h^{-1}$ and since all cyclic subgroups
of order $11$ in $H_0$ are conjugate inside $H_0$ we can write
$(g_0) = hh'(g_0)h'^{-1}h^{-1}$ for some $h'\in H_0$. This implies
$hh'\in N_G((g_0))$. Since $\#N_G((g_0))=\#N_{H_0}((g_0))=55$ (see
the proof of Lemma \ref{orders}), we obtain that
$N_G((g_0))=N_{H_0}((g_0))\subset H_0$, hence $h\in H_0$ and
$h(S_0) = S_0$.  This proves the claim and shows that $H_0$ has
orbits on the set $G(S_0)\setminus \{S_0\}$ of cardinality
divisible by 11. Write
$$G([S_0])=\{[S_0],\, [S_1],\, [S_2],\, ...,\, [S_{r-1}]\}.$$
It is clear that the divisor
$$\calS=\sum_{j=0}^{r-1}S_j$$ is $G$-invariant.
Let $S_0\cdot F_t = b_i$ for $t\in \calO_i$. Then we have
\begin{equation}\label{comp2}
\calF\cdot \calS=(\sum_{j=0}^{r-1}F_j)\cdot (\sum_{j=0}^{r-1}S_j)=rS_0\cdot\sum_{j=0}^{r-1}F_j
=r(1+\sum_{i=1}^{s}b_i\#\calO_i).
\end{equation}
In either case $G\cong M_{11}$ or $M_{22}$, we know $\mu(G)=3$ and hence the two divisors
$\calF$ and $\calS$ are linearly dependent in $\Pic(X_0)$. This implies
$$\calF^2\calS^2 = (\calF\cdot \calS)^2.$$
Substituting from \eqref{comp1}, \eqref{comp2},  we get
\begin{equation}\label{comp3}
r(\sum_{i=1}^sm_i\#\calO_i)\calS^2
=r^2(1+\sum_{i=1}^{s}b_i\#\calO_i)^2.
\end{equation}
Since $\#\calO_i\equiv 0$ mod 11 for all $i$  and $r\equiv
1$ mod 11, the left hand side $\equiv 0 \mod 11$, but
the right hand side $\equiv 1 \mod 11$, a contradiction.

Assume $G_0 = G$. Then the divisor $\calS=S_0$ is $G$-invariant, and
we have a simpler equality
\begin{equation}\label{comp4}
r(\sum_{i=1}^sm_i\#\calO_i)\calS^2
=(1+\sum_{i=1}^{s}b_i\#\calO_i)^2,
\end{equation}
again a contradiction.
\end{proof}

\begin{remark} In \cite{Ko} Kondo proves that the unique supersingular K3 surface
$X$ with Artin invariant 1 admits symplectic automorphism groups
$G\cong M_{11}$ or $G\cong M_{22}$. It follows from the previous
results that any element $g\in G$ of order 11 leaves invariant an
elliptic pencil without a $g$-invariant section. In fact,
according to his construction of $G$ on $X$, one can show that
$\Pic(X)^g\cong U(11)$, hence a $(g)$-invariant elliptic pencil
has only an 11-section.

It is known that the Brauer group of a supersingular K3 surface is isomorphic to the additive group of the field $k$ \cite{Artin}. It is well-known that the group of torsors of an elliptic fibration with a section is isomorphic to the Brauer group. We do not know which  torsors admit a non-trivial automorphism of order $p$ (maybe all?).  We also do not know whether they define elliptic fibrations on the same  surface $X_0$. Note that the latter could happen only for torsors of order divisible by $p =\textup{char}(k)$.  It would be very interesting to see how the three groups
$L_2(11)$, $M_{11}$ and $M_{22}$ sit inside the infinite group $\Aut(X_0)$.
\end{remark}

\begin{remark} It follows from Lemma \ref{lemma} that our surface $X_0$ admits a non-symplectic automorphism of
order $12$. By Remark \ref{goto},  $X_0$ is supersingular with Artin
invariant $\sigma = 1$. It follows from \cite{Ny} that the maximal
order of a non-symplectic isomorphism of a supersingular surface with
Artin invariant $\sigma$ divides $1+p^\sigma$. Thus $12$ is the
maximum possible order. What is the maximum possible non-symplectic
extension of $M_{11}$ or $M_{22}$?
\end{remark}

\begin{remark} A K3 surface may admit a non-symplectic automorphism of
order 11 over any field of characteristic $0$ or $p \ne 2,3,11$. The
well-known example is the surface
$V(x^2+y^3+z^{11}+w^{66})$ in $\bbP(1,6,22,33)$. It is interesting to
know whether there exists a supersingular K3 surface $X$ which admits
a non-symplectic automorphism of order 11. It follows from \cite{Ny}
that, if $p \ne 2$, then 11 must divide $1+p^\sigma$, where $\sigma$ is
the Artin invariant of $X$.
\end{remark}


\end{document}